\documentclass{article}
\usepackage{amssymb,bbm}

\catcode`\à=\active \defà{\`a} \catcode`\À=\active \defÀ{\`A}
\catcode`\á=\active \defá{\'a} \catcode`\Á=\active \defÁ{\'A}
\catcode`\ä=\active \defä{\"a} \catcode`\Ä=\active \defÄ{\"A}
\catcode`\â=\active \defâ{\^a} \catcode`\Â=\active \defÂ{\^A}
\catcode`\å=\active \defå{{\aa}} \catcode`\Å=\active \defÅ{{\AA}}
\catcode`\ç=\active \defç{\c{c}} \catcode`\Ç=\active \defÇ{\c{C}}
\catcode`\è=\active \defè{\`e} \catcode`\È=\active \defÈ{\`E}
\catcode`\é=\active \defé{\'e} \catcode`\É=\active \defÉ{\'E}
\catcode`\ë=\active \defë{\"e} \catcode`\Ë=\active \defË{\"E}
\catcode`\ê=\active \defê{\^e} \catcode`\Ê=\active \defÊ{\^E}
\catcode`\ì=\active \defì{\`{\i}} \catcode`\Ì=\active \defÌ{\`{\I}}
\catcode`\í=\active \defí{\'{\i}} \catcode`\Í=\active \defÍ{\'{\I}}
\catcode`\ï=\active \defï{\"{\i}} \catcode`\Ï=\active \defÏ{\"{\I}}
\catcode`\î=\active \defî{\^{\i}} \catcode`\Î=\active \defÎ{\^{\I}}
\catcode`\ò=\active \defò{\`o} \catcode`\Ò=\active \defÒ{\`O}
\catcode`\ó=\active \defó{\'o} \catcode`\Ó=\active \defÓ{\'O}
\catcode`\ö=\active \defö{\"o} \catcode`\Ö=\active \defÖ{\"O}
\catcode`\ô=\active \defô{\^o} \catcode`\Ô=\active \defÔ{\^O}
\catcode`\ù=\active \defù{\`u} \catcode`\Ù=\active \defÙ{\`U}
\catcode`\ú=\active \defú{\'u} \catcode`\Ú=\active \defÚ{\'U}
\catcode`\ü=\active \defü{\"u} \catcode`\Ü=\active \defÜ{\"U}
\catcode`\û=\active \defû{\^u} \catcode`\Û=\active \defÛ{\^U}
\catcode`\ý=\active \defý{\'y} \catcode`\Ý=\active \defÝ{\'Y}
\catcode`\ÿ=\active \defÿ{\"y} \catcode`\˜=\active \def˜{\"Y}
\catcode`\½=\active \def½{!`}
\catcode`\¾=\active \def¾{?`}
\catcode`\ß=\active \defß{{\ss}}
\newtheorem{lemma}{Lemma}[section]
\newtheorem{theorem}[lemma]{Theorem}
\newtheorem{defi}[lemma]{Definition}
\newtheorem{proposition}[lemma]{Proposition}
\newtheorem{coro}[lemma]{Corollary}
\newenvironment{proof}{
  \noindent\textit{Proof}\ }{\hspace*{\fill}
  \begin{math}\Box\end{math}\medskip}
\newenvironment{proof*}[1]{
  \noindent\textbf{#1\ }}{\hspace*{\fill}
  \begin{math}\Box\end{math}\medskip}
\newcommand{\length}{{\rm length }}

\begin{document}

\title{Straight thin combings} \author{François Dahmani\footnote{The author acknowledges the support of FIM, ETH Zurich.}} \date{} \maketitle

\begin{abstract}
  We introduce a class of thin combings and prove that it has a nice behaviour
  in free constructions.  Among corollaries, we deduce the relative
  metabolicity   
  of Sela's limit groups. 
\end{abstract}

S. Gersten defined the property of metabolicity for finitely presented groups (or even, contractible spaces with bounded geometry in dimension 2) as an attempt of an algebraic approach of word-hyperbolic groups \cite{G1} \cite{G2}. A group $\Gamma$ is said metabolic if its second $\ell_{\infty}$-cohomology group, $H^2_{(\infty)}(\Gamma,A)$ is trivial, for any normed abelian group $A$ (see \cite{G1} for definitions). This condition is stronger than hyperbolicity, but the converse is an open question for groups (for spaces, there are counterexamples exposed in \cite{G2}). 

 In general it is difficult to prove that a given group is metabolic. Free and virtually free groups are the most obvious examples. A consequence of a delicate construction of Rips and Sela is that small cancellation $C'(1/8)$ groups are metabolic. A theorem of Gersten \cite{G2} shows that all hyperbolic surface groups are metabolic, by proving that certain amalgamated free products, over so-called amiable subgroups, preserve the metabolicity. But the problem of deciding if a given subgroup is amiable is neither easy. It is not known in general whether the maximal cyclic subgroups of a metabolic group are amiable. It is neither known whether all the cocompact Kleinian groups are metabolic.

  There is a combinatorial characterisation of metabolicity (for a contractible complex with bounded geometry in dimension 2) given by the existence of a thin combing. A thin combing in the 1-skeleton of a  2-complex is the assignation for each vertex $v$ of a path from the base point to $v$, such that whenever two vertices are neighbors, the combinatorial area of the loop defined by their combings, and by the edge between them is bounded above by an universal constant.

In \cite{G2}, Gersten introduces the relative notion for metabolicity. This amounts to consider the metabolicity of a relative Rips complex (see \cite{Dplms} and below) of a group given with a family of subgroups.   

We prove here that certain relatively hyperbolic groups are relatively metabolic. Our main interest is in Sela's limit groups, and in groups acting freely on $\mathbb{R}^n$-trees. 

\begin{theorem}
Limit groups, and groups that act freely on $\mathbb{R}^n$-trees are metabolic relative to their maximal non cyclic abelian subgroups. 
 
\end{theorem}

This gives an extension to our result in \cite{Dgt}.

We will use the ``d\'evissage'' of such groups (see \cite{VG}), that is, the fact that they split as a graph of groups whose edge groups are cyclic or trivial, and whose vertex groups act freely on $\mathbb{R}^{n-1}$-trees. By performing such splittings, one is left with abelian ``bricks'', after $n+1$ steps.

 Although thin combings have no reason, in general, to behave quietly when performing amalgamations or HNN extensions over a cyclic group, we prove that, if the factors admit so-called \emph{almost prefix-closed, quasi-geodesic} thin combings (what we will call \emph{straight} thin combings), relatively to the cyclic group over which the amalgamation (or HNN) is performed,  then the resulting group admits a straight thin combing, relatively to the cyclic group. This is Theorem \ref{theo;free}, and it occupies the third part of the paper. Moreover, we will prove, in part 2,  that it is possible to expand/collapse straight thin combings when removing/adding a maximal cyclic group in the familly of parabolic subgroups. For example, if a group is metabolic, with a straight thin combing, then it is relatively metabolic, relative to a maximal cyclic subgroup, with a straight thin combing in a suitable conned-off graph, and conversely.  

This allows us to show the existence of a straight thin combing on the Cayley graph of hyperbolic limit groups, and more generally on a suitable conned-off graph of any limit group, where the peripheral subgroups are the maximal abelian non-cyclic subgroups. This gives new examples of metabolic, and of relatively metabolic groups.

 This work followed from a  discussion with A.Szczepa\'nski.

\section{Definitions, thin combings.}

 \subsection{Thin combings}

In the following, by ``path $c$ in a graph $X$'', we mean a continuous map $c$ of a segment $[0, T]$ of $\mathbb{R}$ in $X$, which is simplicial regarding the structure of graph of $\mathbb{R}$ in which vertices are natural numbers. We will say that the path $c:[0,T] \to X$ \emph{starts} at $c(0)$ and \emph{ends} at $c(T)$, where $T$ is the upper bound of the segment on which it is defined. 

If $c$ is a simple path in a graph $X$, from a point $x$ to a point $y$, that is : $c:[0, T] \to X$, $c(0)=x$ and $c(T)=y$, and if a point $z$ is on $c$, that is if there exists a (unique) number $t\in [0,T]$ such that $c(t)=z$, we will denote by $c|_{[x,z]}$ the path $c:[0,t] \to X$ starting at $x$ and ending at $z$. With this vocabulary we will be able to omit the notations of the parameterizations of paths, thus lightening the exposition.

For $\lambda \geq 1$, a $\lambda$-quasigeodesic in a metric space $X$ is a path $c:[0, T] \to X$ that satisfies : $\forall t, t', \frac{|t-t'|}{\lambda} \leq  d_X(c(t),c(t'))$, where $d_X$ represents the distance on $X$. As we choose the path to be simplicial, it is obvious that one always has $d_X(c(t),c(t')) \leq \lambda |t-t'|$, if $\lambda \geq 1$.

If a graph $X$  
is the 1-skeleton of  a simply connected simplicial complex $P$, we define the area of a loop $l$ (generally a triangle or a digone) of $X$ measured in $P$, to be the combinatorial area of the image of $l$ in $P$. If there is no ambiguity on the polyhedron $P$, we write this number $Area(l)$.  We will use only 
polyhedra that are \emph{locally finite around edges}, that is where every edge is contained in finitely many simplicies.

\begin{defi} \label{def;thin}
 
  Let $\mathcal{K}$ be a connected graph, and $v_0$ be a vertex in $\mathcal{K}$. For every number $\lambda$,  
  we denote by $\mathcal{Q}(\lambda,v_0)$ the set of the $\lambda$-quasi-geodesic segments of $\mathcal{K}$ starting at $v_0$. 
 
  A $\lambda$-straight thin combing of $\mathcal{K}$ from $v_0$ is a map $\rho : \mathcal{K}^0 \to \mathcal{Q}(\lambda,v_0)$  
  that assigns to every vertex $v$ of $\mathcal{K}$, a $\lambda$-quasi-geodesic $\rho(v)$,  
  satisfying the two following properties. 
 
  (i) The map $\rho$ is a thin combing, that is: for every $v$, the end point of $\rho(v)$ is $v$, and there is a complex $P$ locally finite around edges, containing the graph $\mathcal{K}$ as its 1-skeleton, there exists $M$ a number, such that for every edge $(v,v')$ of $\mathcal{K}$, the combinatorial area of the triangle defined by the  
  paths $\rho(v)$, $\rho(v')$ and $(v,v')$ is bounded above by $M$, when measured in $P$.  
  We write it $Area (\rho(v),\rho(v'), (v,v') )\leq M$.
 
  (ii) The map $\rho$ is almost prefix-closed, that is, there exists $M>0$ such that for every vertex $v$,  
  and every vertex $v' \in \rho(v)$, the combinatorial area of the digone defined by the paths $\rho(v')$  
  and $\rho(v)|_{[v_0,v']}$ is bounded above by $M$.

  If a group $\Gamma$ acts on the graph $\mathcal{K}$, we say that $\rho$ is a $\Gamma$-thin combing if the areas are measured in a complex $P$ on which the action of $\Gamma$ extends.
 
\end{defi}

 We now establish or recall a few useful lemmas.

\begin{lemma}\label{lem;ext_qg}
  Let $\lambda \geq 1$ be a number. Let $c$ be a $\lambda$-quasigeodesic of a graph, 
  and let $v$ one of its ends. Let $e=(v,v')$ be an edge of the graph
  adjacent to $v$ 
  and such that its other vertex $v'$
  does not lie on $c$. Then the concatenation of $\rho$ and $e$ is a $3 \lambda$-quasigeodesic. 
\end{lemma}

\begin{proof}
  Let $w$ be a vertex of $c$, let $\Delta$ be the distance between $w$ and $v'$ and let $d$ be 
  the length of $c$ between $w$ and $v$. One knows that $d\leq \lambda \times (\Delta+1)$ 
  (quasigeodesy of $c$). As $\lambda \geq 1$ and $\Delta \geq 1$, 
  $d+1\leq 2\lambda\Delta+ \lambda \leq 3\lambda\Delta$, which is what we wanted.   
\end{proof}


\begin{lemma}(Independance of the base point)\label{lem;pt_base}

  Let $\mathcal{K}$ be a graph with an isoperimetric inequality 
  $Area(l)\leq A(\length(l))$ for every simple loop $l$, with areas measured in a complex $P$. Let $v_0$ be 
  in $\mathcal{K}$ and $v_0'$ be another vertex at distance $1$ from $v_0$. 
  We assume that there exists $\rho$, a $\lambda$-straight 
  thin 
  combing of $\mathcal{K}$ from $v_0$, with areas measured in the same complex $P$. Then there exists a $3\lambda$-straigth  
  thin combing of $\mathcal{K}$ from $v_0'$, where areas are measured in the same complex $P$.
 
\end{lemma}

\begin{proof}
  Let us define $\rho'$ as follows. Let $v$ be a vertex of $\mathcal{K}$. 
  If $\rho(v)$ does not contain $v_0'$, then $\rho'(v)$ is the concatenation 
  of the edge $(v_0',v_0)$ in $\mathcal{K}$, with the path $\rho(v)$. If  $\rho(v)$ 
  contains $v_0'$, then $\rho'(v)=\rho(v)|_{[v_o',v]}$. Note that it is well defined 
  since $\rho(v)$ passes only once at $v_0'$, as it is a quasigeodesic.

  Let us check that $\rho'$ is a thin combing. Let $v_1$ and $v_2$ be neighbors in $\mathcal{K}$.
  We can distinguish three cases. 
  First, if $v_0'$ is neither in $\rho(v_1)$ nor in $\rho(v_2)$, then 
  the area of the triangle $( \rho'(v_1), \rho'(v_2), (v_1,v_2))$ is obviously 
  the same as the area of$( \rho(v_1), \rho(v_2), (v_1,v_2))$. 
  Second case, if $v_0'$ is in $\rho(v_1)$ and not in $\rho(v_2)$. 
  By property of quasigeodesy for $\rho(v_1)$, the length of $\rho(v_1)|_{[v_0,v_0']}$ is 
  at most $\lambda$. Consider $z$ the first common point of $\rho(v_1)$ and $\rho(v_2)$ 
  after $v_0'$ on $\rho(v_1)$. 
  By property of thin combing for $\rho$, $\length(\rho(v_1)|_{[v_0', z]}) \leq M$.  Hence, one has  
  $\length(\rho(v_1)|_{[v_0,z]}) \leq \lambda+M$, and $\length(\rho'(v_1)|_{[v_0',z]}) \leq M$.
  By quasi-geodesy of $\rho(v_1)$, the distance between $v_0$ and $z$ is at 
  most $\lambda \times(\lambda+M)$, and therefore,  
  $\length(\rho(v_2)|_{[v_0,z]}) \leq \lambda^2 \times (\lambda +M) $.  
  Therefore the length of the concatenation of 
  $\rho'(v_2)|_{[v_0', z]}$ with $\rho'(v_1)|_{[v_0', z]}$ is at most $\lambda^2 (\lambda+M) +M$.  
  We deduce that $Area(\rho'(v_1), (v_1,v_2),  \rho'(v_2))$ is at most 
  $A(\lambda^2 (\lambda+M) +M) +  Area(\rho'(v_1)|_{[z, v_1]}, (v_1,v_2), \rho'(v_2))|_{[z, v_2]} )$, that 
  is at most $A(\lambda^2 (\lambda+M) +M)+M$. 
  We now turn to the third case when both  $\rho(v_1)$ and $\rho(v_2)$ contain $v_0'$. In this case the triangle 
  $(\rho(v_1), (v_1,v_2),  \rho(v_2))$ is the disjoint union of a digone from $v_0$ to $v_0'$, 
  namely $(\rho(v_1)|_{[v_0,v_0']},  \rho(v_2)|_{[v_0,v_0']} )$,  and a triangle that is nothing else than 
  $(\rho'(v_1), (v_1,v_2),  \rho'(v_2))$, by definition of $\rho'$. Therefore the area of 
  $(\rho'(v_1), (v_1,v_2),  \rho'(v_2))$ is bounded by $M$.

  We now have to check that the thin combing is almost prefix closed. 
  Let $v$ be a vertex, and let $w$ be on $\rho'(v)$. Note that either $w$ is $v_0'$, and there is nothing to prove, 
  or $w$ is on $\rho(v)$.
  First, if $\rho(v)$ does not contain $v_0'$, and $\rho(w)$ does not contain $v_0'$,  
  then the area of the digone  $(\rho'(v)|_{[v_0',w]}, \rho'(w))$ is exactly equal to the area of   
  $(\rho(v)_{[v_0,w]}, \rho(w))$, hence is at most $M$, the constant for the prefix-closedness of $\rho$.
  If $\rho(v)$ does not contain $v_0'$, and $\rho(w)$ contains $v_0'$, then let $w^-$ 
  be the first vertex of $\rho(w)$ after $v_0'$, which is on $\rho(v)$. By assumption that $\rho$ is prefix-closed, 
  the length of $\rho(w)$ between $v_0'$ and $w^-$ is at most $M$.
  Therefore, the length of the loop 
  $(\rho(w)|_{[v_0',w^-]}, (v_0^-, v_0) , \rho(v)|_{[v_0,w^-]})$ is at most $M+1+(M+1)=2M+2$. 
  Hence, the area of the digone $(\rho'(v)|_{[v_0',w]}, \rho'(w))$ is at most $M+A(2M+2)$, 
  where $A$ is the isoperimetric function for $\mathcal{K}$ of the statement.
  If  $\rho(v)$ contains $v_0'$, and  $\rho(w)$ does not, the same argument concludes, with inverted notations.
  Now assume that $\rho(v)$ contain $v_0'$, and  $\rho(w)$ also. 
  Then  the area of the digone  $(\rho'(v)|_{[v_0',w]}, \rho'(w))$ is bounded above by  the area of   
  $(\rho(v)_{[v_0,w]}, \rho(w))$, and is at most $M$.

  Finally, we have to check that for every $v$, $\rho'(v)$ is a $3\lambda$-quasigeodesic. 
  If $v_0'$ is in $\rho(v)$, it is obvious because $\rho'(v)$ is a subpath of $\rho(v)$. 
  If $v_0'$ is not in $\rho(v)$, this is a consequence of Lemma \ref{lem;ext_qg}.

\end{proof}

\subsection{Hyperbolicity, fineness, and $\Gamma$-graphs of finite type}

 Following B.Bowditch in \cite{Brel}, we define the notion of \emph{$\Gamma$-graph of finite type}, for a group $\Gamma$, to be a 
 graph on which $\Gamma$ acts simplicially, with finitely many orbits of edges, and such that the stabilizer of every edge is finite.
  In such a case, we call the stabilisers of the vertices of infinite valence, the \emph{peripheral subgroups} of 
 $\Gamma$ acting on the graph.

 We recall also the definition of \emph{fineness} of a graph (also from \cite{Brel}). A graph is said \emph{fine} if for every edge $e$ and every number $L$, the set of simple loops containing $e$ and of length $L$ is finite. 

 By definition (again, found in  \cite{Brel}), for a group $\Gamma$, the existence of a hyperbolic fine $\Gamma$-graph of finite type implies the relative hyperbolicity of $\Gamma$ (relative to the stabilisers of vertices of infinite valence). To read more about this, the reader is of course referred to \cite{Brel}, and also to \cite{Dt}. We will make use of  
the existence of a Relative Rips complex for relatively hyperbolic groups, proven in \cite{Dplms}.

\begin{proposition} (Theorem 6.2 in  \cite{Dplms})

Let $\Gamma$ be a group, and $\mathcal{K}$ be a $\Gamma$-graph of finite type, hyperbolic and fine. Then, there exists a simplicial complex $P(\mathcal{K})$, with a simplicial action of $\Gamma$ of finite dimension, locally finite except at the vertices, and contractible, whose 1-skeleton is a  $\Gamma$-graph of finite type containing $\mathcal{K}$. 

\end{proposition} 

In  \cite{Dplms}, such a complex is constructed, and is called a relative Rips complex.
 
\begin{coro}

 For any relative Rips complex, there is an isoperimetric function, that is a function $f$ such that any loop of length $R$ bounds an area of at most $f(R)$.

\end{coro}

It is a direct consequence of the finiteness of simple loops of given length up to translation, and that the complex is simply connected. Note that we only use the fineness. With the hyperbolicity, it is in fact easy to see that this function can be taken linear.

\subsection{Farb's conning off}


  Let $\Gamma$ be a group, and 
   $\mathcal{K}$ a $\Gamma$-graph of finite type. 
  Let $H$ be a subgroup of $\Gamma$. 
  We choose $F$ a family of representatives of the $\Gamma$-orbits of vertices of $\mathcal{K}$. 
  For every left coset $\gamma H$ of $H$, we add a vertex $\tilde{v}_{\gamma H}$, 
  and edges $(\tilde{v}_{\gamma H},\gamma h v)$ for every $h\in H$, and every $v\in F$.

  Following Farb in \cite{F}, we call the resulting graph the \emph{coned-off} graph 
  of $\mathcal{K}$ along along the orbits (on the right) of $H$. 
  When there is no ambiguity we write $\widehat{\mathcal{K}}$ for this graph. We call the additionnal vertices of 
  $\widehat{\mathcal{K}} \setminus \mathcal{K}$, the coset-vertices. 

  If one has a polyhedron $P$ containing $\mathcal{K}$ in its 1-skeleton, then we define $\hat{P}$ to be the 
  complex whose vertices are those 
  of $\mathcal{K}$, and whose vertices are those of $P$ and the coset-vertices, whose edges are the edges of $\widehat{\mathcal{K}}$, 
   whose 2-simplices are those of $P$ and the additional 2-simplices 
  containing a coset-vertex of $\mathcal{K}$ and three edges of $\widehat{P}$, and whose higher dimensional 
  simplices are those of $P$.  Note that if $P$ is locally finite around the edges, then $\hat{P}$ also.

\subsection{Property of BCP}

 When looking at the coset-vertices of a conned-off graph, the property of fineness implies the property of BCP that we now explain.

Consider a simple path $\hat{c}:[0,T]\to \widehat{\mathcal{K}}$,  in $\widehat{\mathcal{K}}$, and 
  $\tilde{v}=\hat{c}(t)$ a coset-vertex on it. The vertex $\hat{c}(t-1)$ is a vertex of 
  $\mathcal{K}\subset \widehat{\mathcal{K}}$, and we call it the entering vertex 
  of $\hat{c}$ in the coset of $\tilde{v}$. Similarily, we call $\hat{c}(t+1)$ the exiting 
  vertex of $\hat{c}$ in the coset of $\tilde{v}$.


  We say that the pair $(\Gamma, H)$ satisfies the property of 
  \emph{Bounded Coset Penetration} (property of BCP) in $\mathcal{K}$ if, 
  for all number $\lambda$, there exists a number $r_{BCP}(\lambda)$ such that 
  whenever two $\lambda$-quasigeodesics $\hat{c}$  
  and $\hat{c}'$ of $\widehat{\mathcal{K}}$ 
  with same end points pass on a same coset-vertex $\tilde{v}$, the entering 
  vertices of both are at distance at most $r_{BCP}(\lambda)$ in 
  the metric of $\mathcal{K}$.

  One has the proposition (we refer to \cite{Dplms} or \cite{Dt}-Appendix for a proof).
 
\begin{proposition}

 If a conned-off graph is hyperbolic and fine, then it satisfies the property of BCP (as defined below).

\end{proposition} 
 






\subsection{Reductions, and extensions of paths}

  We now define the \emph{Farb's reduction} of a path of $\mathcal{K}$. 
  Given a path $c: [t_0,T] \to \mathcal{K}$ in $\mathcal{K}$, parameterised by arc-length,
  we look for a vertex $c(t)$, $t\geq t_0+2$ such that $c(t)$ is a neighbor of a same vertex-coset than $c(0)$. 
  If such a number $t$  exists, we choose the last one among them, $t_1$, and we replace $c|_{[t_0,t_1]}$ by the 
  path of two edges joining $c(t_0)$ and $c(t_1)$ 
  passing at the vertex-coset. If there is no such number $t$, we simply note $t_1=t_0+1$. Then we consider $c|_{t_1,T}$ 
  and we do the same, hence defing $t_2\geq t_1 +1$. Eventually we arrive at the end of the path, that is $t_k =T$ for some $k$. 
  Then the obtained path, with the consecutive replacements, is the Farb reduction of $c$, and we note it $\hat{c}$.


  Conversely, we define the \emph{extension} of a simple  
  path of $\widehat{\mathcal{K}}$. 
  Let $c$ be a simple path in $\widehat{\mathcal{K}}$, between two vertices of $\mathcal{K}$, noted $v$ and $v'$.
  Let $w_1,\dots, w_k$ be the cosets vertices lying on $c$ in this order from $v$ to $v'$. 
  The edges of $c$ are all in $\mathcal{K}$ except precisely $2k$ edges, 
  namely those that are before or after a vertex $w_i$. For each $i$ we note $w_i^-$, 
  respectively $w_i^+$, the vertices of $c$ lying just before, respectively 
  just after $w_i$. Note that these vertices are in 
  $\mathcal{K}$. Note also that there exists a conjugate of an element of $H$, 
  $h^\gamma$ such that $h_i^{\gamma_i} w_i^- = w_i^+$. 
  We choose in $\mathcal{K}$ the path $c_i$ between $w_i^-$ and $w_i^+$ corresponding to the 
  path given by the thin combing of $H$ for $h$ in the graph $\gamma_i Cay(H) \subset \mathcal{K}$. 
  The path $\check{c}$ consists of the concatenation: 
  $\check{c} = c|_{[v,w_1^-]}\, c_1 \,c|_{[w_1^+, w_2^-]}\, c_2 \, \dots c_{k-1}\,  c|_{[w_{k-1}^+, w_k^-]}\, c_k \, c|_{[w_{k}^+, v']}$.

%
%

\subsection{Lemmas of quasi-geodesy}

We say that a subgroup of a group $\Gamma$ acting on a graph $X$ is $X$-quasi-convex in $X$ if the orbit of a certain point in $X$ is quasi-convex in $X$.   

 We will make use of two useful lemmas. Independantly, and simultaneously, C.~Dru\c{t}u and M.~Sapir proved a stronger version in \cite{DS} (Theorem 1.12). We give our proof below so that the exposition of this important part is self-contained.

\begin{lemma}\label{lem;qg1}
  Let $\Gamma$ be a group, and  $\mathcal{K}$ be a hyperbolic and fine $\Gamma$-graph of finite type.   
  With the notations above, let $H$ be a $\mathcal{K}$-quasi-convex  subgroup of $\Gamma$, and $\widehat{\mathcal{K}}$ be the conned-off graph of $\mathcal{K}$ along the orbits of $H$. 

  If $\widehat{\mathcal{K}}$ satisfies the property of BCP,  
  then for all $\lambda$ 
  there exists $\mu$ such that for all $\lambda$-quasigeodesic 
  $c$ of $\mathcal{K}$, the Farb reduction $\hat{c}$ is a $\mu$-quasigeodesic 
  of $\widehat{\mathcal{K}}$.
\end{lemma}

\begin{proof}

\emph{Remark.}  The hypothesis of quasi-convexity is always fulfiled if $H$ is cyclic.  In fact, if $H$ is parabolic, that is  a subgroup of an element of the family $\mathcal{G}$, then there is nothing to prove. If $\mathcal{H}$ is not parabolic, then a generator is either a eliptic isometry or an hyperbolic isometry of $\mathcal{K}$. If it is elliptic, it has finite order, and $H$ is finite, again the result is obvious in this case. If it is hyperbolic, it stabilises a bi-geodesic, and it follows that the orbit of a base point of $\mathcal{K}$ under $H$ is  quasi-convex in $\mathcal{K}$. Note also that in this case, if $H$ is maximal cyclic, then $\widehat{\mathcal{K}}$ satisfies the property of BCP, by a proposition of \cite{Dgt} (Lemma 4.4).

Let $C$ be a constance of quasi-convexity for $H$ in $\mathcal{K}$. 
 
  Let $v$ and $v'$ be two points $\hat{c}$, not coset-vertex, and let $g$ be a geodesic of $\widehat{\mathcal{K}}$ 
  from $v$ to $v'$. 
  Let $\Delta$ be the length of $g$, and $d$ be the length of $\hat{c}$ between 
  $v$ and $v'$. We have to find an upper bound for $d$ in terms of $\Delta$. For that, it is enough, up to 
  subdivision of paths, to restrict to the case where $g$ and $\hat{c}$ do not travel in a same vertex-coset.

  Let us consider $\check{g}$ the extension of $g$, in $\mathcal{K}$. If the longest travel of $g$ in vertex-cosets 
  has length $T(g)$, then $\length(\check{g}) \leq \Delta + \Delta \times T(g)$, 
  and therefore the distance from $v$ to $v'$ in $\mathcal{K}$ is at most  $\Delta + \Delta \times T(g)$, and 
  by quasi-geodesy of $c$, $\length(c|_{[v,v']}) \leq  \lambda (\Delta + \Delta \times T(g))$. Therefore, 
  as reductions contract distances, $d \leq  \lambda (\Delta + \Delta \times T(g))$. 
  It is enough, in order to conclude,  
  to give an upper bound for  $T(g)$ that  depends only on $\Delta$. 

  The path $\check{g}$ consists of the concatenation of at most $2\Delta$ geodesic segment and  
  $C$-quasi-geodesic segments, where $C$ is the constant of quasi-convexity of the subgroup $H$ in $\mathcal{K}$. 
  Let $g_0$ be a longest travel in a $H$-coset in $\check{g}$, 
  the distance between its end points of $\mathcal{K}$ is $T(g)$. Let $g_1$ be another travel 
  in $H$-coset of $g$, different from $g_0$. For all $L$, the segment $g_1$ has no subsegment 
  of length $C\times (r_{BCP}(1) +2L)$ 
  that stays within distance $L$ from $g_0$, because if it had one, we would concludes a contradiction for the property 
  BCP of $\widehat{\mathcal{K}}$, as there would be segments between two coset-vertices of length at 
  most $L$ with exiting vertices at distance at least $r_{BCP}(1)$. 
  Therefore, by hyperbolicity for the polygon defined by the concatenations of 
  the subsegments of $\check{g}$ and of $c$, there is a constant $R$ depending only on $\Delta$, $C$, 
  $\lambda$ and $\delta$, the hyperbolicity constant 
  of $\mathcal{K}$, such that  
  for all $n$, there exists $N(n)$ such that if $T(g) \geq N(n)$, the paths $g_0$ and $c$ remain $R$-close to 
  each other for a length at least  $n$. Let us consider the first point of $c$ of this subsegment, to be $c(t_0)$. 
  Let us consider the parameterisation of paths:  $c:[t_0,t_0+n] \to \mathcal{K}$, and $g_0:[0,n] \to \mathcal{K}$.
  Now by the lemma of conical stability of quasi-geodesics, see Proposition 1.11 in \cite{Di}, 
  there is a constant $S$ depending only on $R$, $C$, $\lambda$, and $\mathcal{K}$, such that the family of segments $[g_0(k),c(t_0+k)]$, $k=1\dots n$ contains at most $S$ elemtents modulo left-translation. 

Therefore, for all $n'$ there exists $n$ such that if $T(g) \geq N(n)$,  the family of segments $[g_0(k),c(t_0+k)]$, $k=1 \dots n$ contains $n'$ times the same segment up to translation by elements of $\Gamma$. Let $[g_0(k_i),c(t_0+k_i)]$, $i=1 \dots  n'$ be these segments. By assumption of $g_0$, there are elements of a conjugate of $H$ in $\Gamma$, say, $h_2 \dots h_{n'}$, such that $g_0(k_i)=h_i g_0(k_1)$, for $i\geq 2$. Therefore, for all $i \geq 2$, $c(t_0+k_i)=h_i c(t_0)$. 
  We see that,
  if $T(g)$ is greater than a certain constant (so that, the number $n'$ above is greater than $2r_{BCP}(C)$), depending only on $\Delta$, $C$, $\lambda$ and $\mathcal{K}$, the path $\hat{c}$ enters and exists the same vertex-coset than the vertex coset associated to $g_0$. This was excluded.

\end{proof}

\begin{lemma}\label{lem;qg2}

  With the same assumptions as in the previous lemma, for all $\mu$, 
  there exists $\nu$ such that for all $\mu$-quasigeodesic 
  $\hat{c}$ of $\widehat{\mathcal{K}}$, its extension $\check{c}$ is a $\nu$-quasigeodesic 
  of $\mathcal{K}$.
\end{lemma}
   
\begin{proof}

  Let us consider $x$ and $y$ on $\check{c}$ at distance $\Delta$ in  
  $\mathcal{K}$, and denote by $d$ the length of $\check{c}$ between them. 
  Let $g$ be a geodesic between $x$ and $y$ in $\mathcal{K}$. By the previous lemma, there 
  is a constant $\lambda$ depending only on the spaces such that $\hat{g}$ 
  is a $\lambda$-quasigeodesic of $\widehat{\mathcal{K}}$. Let $\eta = \max\{\lambda,\mu\}$. 
  By the property of BCP for $\hat{c}$ and $\hat{g}$, the path $\hat{c}|_{[x,y]}$ makes no travel in a coset-vertex that is 
  longer than $r_{bcp}(\eta)$ than the travel of $\hat{g}$ in the same coset-vertex.  
  Let us give a upper bound for $d=\length (\check{c}|_{[x,y]})$. It is less than the length  of  
  $c|_{[x,y]}$ plus the sum of the lengths of the travels of   $c|_{[x,y]}$ in coset vertices. Therefore, $d$, 
  is at most $\length (\hat{c}|_{[x,y]})$ plus the sum of the travels of $\hat{g}$ 
  ( let $T$ be this number) plus  $m \times (2r_{BCP}(\eta)$, where $m$ 
  is the number of coset-vertices encontered by $\hat{g}$. In particular $m$ is less than the length of $\hat{g}$.  
  We have $d\leq \length(\hat{c}|_{[x,y]}) +T + 2r_{BCP}(\eta) \length(\hat{g})$. 
  Now notice that, by the definitions of quasigeodesics, 
  $  \length(\hat{c}|_{[x,y]}) \leq \frac{\mu}{\lambda} \length(\hat{g})$. Therefore,  
  $d\leq (\frac{\mu}{\lambda} +2r_{BCP}(\eta) ) \times  \length(\hat{g}) + T$. 
  As $\frac{\length(\hat{g})}{3} \leq \length(g)=\Delta$,  
  $d\leq 3\times (\frac{\mu}{\lambda} +2r_{BCP}(\eta)) \Delta$. This proves the lemma.
\end{proof}

\section{Equivalence Lemmas, reductions and extensions of straight thin combings.}

%
%
\subsection{Reductions of straight thin combings}

We state and prove the two central Lemmas of this paper.

\begin{lemma}\label{lem;contrac}
  Let $\Gamma$ be a group, and $X$ be hyperbolic fine $\Gamma$-graph of finite type. 
  Let $P(X)$ be a reltaive Rips polyhedron for $X$, and let $\mathcal{K}$ be its 1-skeleton. 
  Let $v_0 \in \mathcal{K}^0$ be a base point in $\mathcal{K}$.
  Assume that there exists $\lambda$ a number, and $\rho$ is a
  $\Gamma$-$\lambda$-straight  
  thin combing of $\mathcal{K}$ from
  $v_0$, with areas measured in $P$.
  
  Let $H$ be a $\mathcal{K}$-quasi-convex  
subroup of $\Gamma$ and  $\widehat{\mathcal{K}}$ be
  the coned-off of $\mathcal{K}$ along the orbits of $H$ in $\Gamma v_0$.
  Assume that $\mathcal{\widehat{\mathcal{K}}}$ satisfies the property of BCP.

  Then, there exists a number $\mu > 0$ and a $\Gamma$-$\mu$-straight 
  thin combing of $\widehat{\mathcal{K}}$ from $\widehat{v_0}$, with areas measured in $\hat{P}$.   
\end{lemma}

\begin{proof}
  We define $\hat{\rho}$ as follows. If $v$ is in $\mathcal{K^0} \subset
  \mathcal{\widehat{\mathcal{K}^{}}}^0$, $\hat{\rho} ( v ) = \widehat{\rho ( v
  )}$ in the sense of the Farb reduction of paths. If $v$ is in
  $\widehat{\mathcal{K}^0} \mathcal{\setminus \mathcal{K}}^0$, then, let us choose a geodesic $ [v_0, v]$ and let $v'$ be the point before $v$ on it. We choose   $\hat{\rho} ( v ) =  \hat{\rho} (v')$.

  We start by noticing that,  by Lemma \ref{lem;qg1}, there is a
  constant $\mu$ such that, for all vertex $v$
  in $\widehat{\mathcal{K}^0}$, the path $\hat{\rho} ( v )$ ends at $v$,  
  and is 
  a $\mu$-quasi-geodesic. 

  We have to prove that $\hat{\rho}$ is a thin combing,  
  and almoxt prefix-closed. 

  Let $M$ be a constant of Definition \ref{def;thin} for $\rho$. We first prove:

  \begin{lemma}\label{lem;tech1}

    Let $\lambda$ be a number. Then there exists a constant $\chi(\lambda)$ depending only 
   on $M$, $\lambda$, $\mathcal{K}$ and $\widehat{\mathcal{K}}$ satisfying the following. 

    For all pairs $(c,c')$ of $\lambda$-quasi-geodesics of $\mathcal{K}$,  with same end points 
    such that $Area(c, c')\leq M$, 
    one has  $Area(\hat{c}, \hat{c'}) \leq \chi(\lambda)$, 
    where $\hat{c}$ and $\hat{c'}$ are the Farb reductions of paths.
  \end{lemma}
  
  \begin{proof*}[Proof.]

    The graph $\mathcal{K}$ is included in the graph $\widehat{\mathcal{K}}$, 
    therefore, the path  $c$ is a simple path of   $\widehat{\mathcal{K}}$.  
    As in $\mathcal{K}$, the area  
    $Area_{\mathcal{K}} (c,c')$ is at most $M$, the inclusion of graphs guaranties that the area in $\widehat{\mathcal{K}}$  
    of this digone is at most $M$: $Area_{\widehat{\mathcal{K}}} (c,c') \leq M$.  Let $\mu$ be a constant 
    of quasi-geodicity given by 
    Lemma \ref{lem;qg1}, for $\lambda$.

    Let $v$ be a vertex coset of $\widehat{\mathcal{K}}$, on $\hat{c}$ but not on $\hat{c'}$. By property of BCP,  
    the path $c$ travels at most for a distance of $r_{BCP}(\mu)$ in this coset-vertex. By quasi-geodesy of $\hat{c}$, this  
    means that the subpath $t$ of $c$ between the entering point and the exiting point in $v$ is at most  
    $\lambda \times r_{BCP}(\mu)$. If we denote by $c_v$ the path of $\widehat{\mathcal{K}}$ obtained from $c$ by the  
    replacement of $t$ by the two edge path passing by $v$, we see that  
    $Area_{\widehat{\mathcal{K}}} (c,c_v) \leq A_{\widehat{\mathcal{K}}} (\lambda r_{BCP}(\mu)+2)$, where  
    $A_{\widehat{\mathcal{K}}}$ is an isoperimetric function for $\widehat{\mathcal{K}}$.  
 
   Therefore $Area_{\widehat{\mathcal{K}}} (c_v,c') \leq M+ A_{\widehat{\mathcal{K}}} (\lambda r_{BCP}(\mu))+2$, 
   and if we denote by  
   $\tilde{c}$ the path obtained by reducing $c$ on every vertex-coset on $\hat{c}$ that are not on $\hat{c'}$ 
   (by assumption on  
   the area betwen $c$ and $c'$ there are at most $M$ of them), we get that  
   $Area_{\widehat{\mathcal{K}}} (\tilde{c},c') \leq M+ M\times A_{\widehat{\mathcal{K}}} (\lambda r_{BCP}(\mu)+2)$.  
 
   In the same way, if $\tilde{c'}$ is the path obtained from $c'$ by reducing it at every coset vertex not on $c$,   
   $Area_{\widehat{\mathcal{K}}} (\tilde{c},\tilde{c'}) 
   \leq M+ 2M\times A_{\widehat{\mathcal{K}}} (\lambda r_{BCP}(\mu)+2)$.   Let $N$ be this upper bound.

    Among the coset-vertices on both $\hat{c}$ and $\hat{c}'$, there are at most $M$ such that  
    that $c$ and $c'$ does not enter or exit the coset-vertex at the same point, by preperty of thin combing for $c$ and $c'$.
     Let $\dot{v}_1 \dots \dot{v}_m$ ($m\leq M$) 
    be these consecutive coset-vertices. For all $i$, there is a path from the entering point of $\tilde{c}$ in $\dot{v_i}$ 
    to the entering point of $\tilde{c}'$, consisting of a subsegment of   $\tilde{c}$  and one of $\tilde{c}'$, and that 
    is not longer than $N$. 
    Indeed if it was not true, the area $Area_{\widehat{\mathcal{K}}} (\tilde{c},\tilde{c}')$ would be greater 
    than $N$, contradicting its definition. A similar statement is true for exiting points. 
    Therefore, if $\tilde{c}_{v_i}$ and $\tilde{c}'_{v_i}$ are the path obtained from 
    $\tilde{c}$ and $\tilde{c}'$ by reducting the travel in the coset of $v_i$, the area 
    $Area_{\widehat{\mathcal{K}}} (\tilde{c}_{v_i},\tilde{c}'_{v_i})$ is at most $N + 2 \times A_{\mathcal{K}} (N+2)$. 

    If $\dot{c}$ and $\dot{c}'$ are the paths obtained from $\tilde{c}$ and $\tilde{c}'$ by 
    reducting the travels in every coset of $v_1\dots v_m$, the area 
    $Area_{\widehat{\mathcal{K}}} (\dot{c}, \dot{c}') \leq N+ 2M\times A_{\mathcal{K}} (N+2) $. 
    Now in every coset-vertex of $\hat{c}$ and $\hat{c}'$ that is not aready on $\dot{c}$ and $\dot{c}'$, 
    the entering points and the exiting points of the two paths are the same. 
    Therefore, the reduction does not change the area, and we conclude that 
    $Area_{\widehat{\mathcal{K}}} (\hat{c}, \hat{c}') \leq N+ 2M\times A_{\mathcal{K}} (N+2) $, that is, 
   less than a constant depending on $M$, $\lambda$, $\mathcal{K}$ and $\widehat{\mathcal{K}}$.

  \end{proof*}  

  To prove that $\hat{\rho}$ is a thin combing, we have to bound the area of the all triangles 
  $(\hat{\rho}(v), e, \hat{\rho}(w))$, where $e=(v,w)$ is an edge, and $v$ and $w$ are 
  vertices of $\widehat{\mathcal{K}}$. We have three case to consider, namely if both $v$ 
  and $w$ are vertices of $\mathcal{K}\subset \widehat{\mathcal{K}}$, when only $v$ is 
  a vertex of $\mathcal{K}$, and when none of them are in $\mathcal{K}$. 
  We can rule out the last case by noticing that the minimal distance between vertices of 
  $\widehat{\mathcal{K}} \setminus \mathcal{K}$ is at least $2$ (obvious by construction).  

  Let us bound the area of a triangle 
  $(\hat{\rho}(v), e, \hat{\rho}(w))$, where $v$ and $w$ are 
  vertices of $\widehat{\mathcal{K}}$. If $(\hat{\rho}(v)$ contains $w$ and conversely, 
  the area is the area of the digone 
  $(\hat{\rho}(v), \hat{\rho}(w)|_{[\widehat{v_0},v]})$,  and, 
  from Lemma \ref{lem;tech1}, it is bounded above by $\chi(\lambda)$.  
  In the other case, up to a change of names, $w$ is not in $\hat{\rho}(v)$. 
  By Lemma \ref{lem;ext_qg}  the concatenation of  $\hat{\rho}(v)$ with $e$ is a 
  $ 3\mu $-quasi-geodesic.  By Lemma \ref{lem;tech1}, the area of the triangle 
  is therefore at most $\chi(3\lambda)$.  

  In the case where $w$ is a coset-vertex, that is a vertex in 
  $\widehat{\mathcal{K}} \setminus \mathcal{K}$ 
  we use the almost prefix closed property for the straight thin combing $\rho$. Let $w'$ be as in the definition 
  of $\hat{\rho(w)}$; we have $ \hat{\rho}(w) =\widehat{\rho(w')}|_{[v_0,w]}$. Let $w^-$ be the vertex of   
  $\hat{\rho}(w)$ lying next to $w$ in  $\hat{\rho}(w)$. Then  $w^-$ is a vertex of $\mathcal{K}$, as no 
  pair of coset-vertices are adjacent. By prefix-closed property for $\rho$, the area of the digone 
  $(\rho(w^-), \rho (w')|_{[v_0,w^-]})$ is at most $M$. Therefore, by Lemma \ref{lem;tech1}, 
  the area of the triangle  defined by  
  $\hat{\rho}(w)$,  $\hat{\rho}(w^-)$ and the edge $(w,w^-)$ is at most   $\chi(3\lambda)$. 
 
  We turn now to $\hat{\rho}(v)$. If this path contains $w$, then, denote by $v^-$ 
  the vertex that come just before $w$ on it, by the same arguments, the area defined by the trianble $\hat{\rho}(v)$,  
  $\hat{\rho}(v^-)$ and the path $\hat{\rho}(v)|_{[v^-,v]}$ is at most   $\chi(3\lambda)$.  
  Moreover, by the property of $BCP$, the vertices $v^-$ and $w^-$ are at distance at 
  most $r_{BCP}(\mu)$ from each other in    $\mathcal{K}$. 
  Let $c$ be a path from $v^-$ to $w^-$ in $\mathcal{K}$ of length less than $r_{BCP}(\mu)$.
  According to the upper bound on areas of triangles whose vertices are in $\mathcal{K}$ we obtained above,  
  and  by triangular inequality for areas, the triangle defined by $\hat{\rho}(w^-)$, $\hat{\rho}(v^-)$, and  $c$ is at most 
  $r_{BCP}(\mu) \times \chi(3\lambda)$. 
  Finally, the area of the loop  
  defined by $c$, $(w^-, w)$ and $(w, v^-)$ is at most $A_{\mathcal{K} } (r_{BCP}(\mu) +2 )$  
  where $A_{\mathcal{K}}$ is an isoperimetric 
  function for  $\mathcal{K}$.
  Therefore the area of the triangle defined by 
  $\hat{\rho}(w), \hat{\rho}(v)$ and the edge $(v,w)$ is at most  
  $r_{BCP}(\mu) \times \chi(3\lambda)                  
  + A_{\mathcal{K} } (r_{BCP}(\mu) +2 )$, 
  which is an universal constant, depending only on $\rho$ and $\widehat{\mathcal{K}}$.

  If now the path $\hat{\rho}(v)$ does not contain $w$, it suffices to perform the same argument 
  with $v^- = v$ to get the majoration of the area.

  This proves that $\hat{\rho}$ is a thin combing. We now turn to the almost prefix-closedness property. 

  Let $v$ be a vertex of $\mathcal{K}$ and $w$ be a vertex of $\hat{\rho}(v)$. We have to bound the area of 
  the digone defined by $\hat{\rho}(w)$ and  $\hat{\rho}(v)|_{[v_0,w]}$.

  If both vertices are in $\mathcal{K}$, then it is a direct corollary of the property of almost prefix-closedness 
  for $\rho$ and of Lemma \ref{lem;tech1}. Let us assume now that $v$ is a coset-vertex, let $v^-$ be the last vertex of 
  $\hat{\rho}(v)$ before $v$. As in the argument above, we know that the area of the digone 
  $(\hat{\rho}(v)|_{[v_0,v^-]}, \hat{\rho}(v^-))$ is uniformly bounded by $\chi(\lambda)$.  

  If $w$ is a vertex on $\hat{\rho}(v)$, and assume it is not a vertex-coset,  then,  the area of the digone 
  $(\hat{\rho}(v^-)|_{[v_0,w]}, \hat{\rho}(w))$ is bounded by  $\chi(\lambda)$.  
  Therefore $Area( \hat{\rho}(v)|_{[v_0,v^-]},\hat{\rho}(w))$ 
  is at most $2\times  \chi(\lambda)$. If now $w$ is a vertex-coset, then let $w^-$ be the last vertex of 
  $\hat{\rho}(w)$ before $w$, and $w'^-$ the last vertex of 
  $\hat{\rho}(v)$ before $w$. By the property of $BCP$ we know that the distance in $\mathcal{K}$ 
  between $w^-$ and $w'^-$ is at 
  most $r_{BCP}(\mu)$, and therefore, there is a path $p$ in $\mathcal{K}$, of length at most $r_{BCP}(\mu)$, and such that  
  $Area(\hat{\rho}(w^- ), \hat{\rho}(w'^-), p) \leq  r_{BCP}(\mu) \times  \chi(3\lambda)$.  Moreover, 
  the area of the loop $((w,w^-) p (w'^-,w)) $ 
  is at most $A_{\widehat{\mathcal{K}}} ( r_{BCP}(\mu) +2)$, where $A_{\widehat{\mathcal{K}}}$ is an isoperimetric function 
  for $\widehat{\mathcal{K}}$. Now we notice that, from the cases treated above, 
  $Area(\hat{\rho}(w^- ), \hat{\rho}(w)|_{[v_0,w^- ]} )\leq  \chi(\lambda)$   and 
  $Area(\hat{\rho}(w'^- ),\hat{\rho}(v)|_{[v_0,w'^- ]}) \leq  \chi(\lambda)$. 
  By triangular inequality for the areas, one concludes: 
  $Area(\hat{\rho}(w), \hat{\rho}(v)|_{[v_0,w]} ) \leq 2\times  \chi(\lambda)  +  r_{BCP}(\mu)\times \chi(3\lambda)  + A_{\hat{\mathcal{K}}} (r_{BCP}(\mu) +2)$.

\end{proof}

%
%

\subsection{Extensions of straight thin combings}

Let us turn now to the converse operation.

\begin{lemma}\label{lem;ext}
  Let $\Gamma$ be a group, and , $X$ be a hyperbolic and fine $\Gamma$-graph of finite type. 
  Let $P(X)$ be a relative Rips polyhedron and $\mathcal{K}$ its 1-skeleton. 
  Let $v_0 \in \mathcal{K}^0$ be a base point in $\mathcal{K}$.
  
  Let $H$ be a quasi-convex 
  subroup of $\Gamma$ admitting a straight thin combing, and  $\widehat{\mathcal{K}}$ be the
  coned-off of $\mathcal{K}$ along the orbits of $H$ in $\Gamma v_0$. We assume
  that $\mathcal{\widehat{\mathcal{K}}}$ satisfies the property of BCP.

  Assume that there exists a number $\mu > 0$ and $\rho$ a $\Gamma$-$\mu$-straight   thin combing of $\widehat{\mathcal{K}}$ from $\widehat{v_0}$, with areas measured in  $\widehat{P(X)}$. 
  Then
  there exists a number $\lambda$, and a $\lambda$-straight 
  thin combing of $\mathcal{K}$ from $v_0$, with areas measured in the subcomplex $P(X)$ of $\widehat{P(X)}$. 
\end{lemma}

\begin{proof}
  
  We define the candidate $\check{\rho}$ as follows. Let $v$ be a vertex of $\mathcal{K}$, It is also in the graph 
  $\widehat{\mathcal{K}}$, and $\check{\rho} (v)$ is  simply $\check{\rho (v)}$, the extension of $\rho (v)$.

  We begin by a lemma.

  \begin{lemma}
    Let $c$ and $c'$ be two simple paths of $\widehat{\mathcal{K}}$ with same end points 
    such that $Area(c, c')\leq M$. Then 
    $Area(\check{c}, \check{c'}) \leq  M\times (A_{\mathcal{K}} (M\times (1+  r_{BCP}(\lambda))) + 2 r_{BCP}(\lambda)  M_H )$, 
    where   $A_{\mathcal{K}}$ is a chosen 
    isoperimetric function of $\mathcal{K}$, and 
    where $\check{c}$ and $\check{c'}$ are the expansions of paths in $\mathcal{K}$.
  \end{lemma}

  \begin{proof*}[Proof.]

    Let $\Delta$ be the set of cosets-vertices of $c$ or $c'$ such that $c$ and $c'$ disagree on the 
    entering point or the exiting point.  
    Note that by assumption on the area of the digone $(c,c')$, 
    the family $\Delta$  has cardinality at most $M$. 
  
    Let $c_1,\dots, c_m$ be the consecutive maximal subsegments that are common to $c$ and $c'$.  
    For all $i$, there is a digone $D_i$ between the end of $c_i$ and the beginning of $c_{i+1}$, consisting 
    of a subsegment of $c$ and a subsegment of $c'$. By assumption, the area of each such digone in $\mathcal{K}$ 
    is at most $M$, the constant of the thin combing, and also, there are at most $M$ such digones, $m\leq M$. 
    By definition, 
    the digone $(\check{c}, \check{c'})$ consists 
    of the union of the $\check(c_i)$ and of the $\check(D_i)$, hence, we only need to bound the area of each $\check(D_i)$.

    If $w$ is an element of $\Delta$, that is, a vertex-coset 
    in one of the loops $D_i$. As $D_i$ is the concatenation of two $\lambda$-quasi-geodesics, by property of BCP, 
    the distance between the two neighbors of $w$ in $D_i$  
    in $\mathcal{K}$ is at most $r_{BCP}(\lambda)$. 

    We see that $\check(D_i)$ consists in a loop of length at most 
    $\length(D_i) + M \times r_{BCP}(\lambda) = M\times (1+  r_{BCP}(\lambda))$ attached to two pairs of 
    paths given by the thin combing of $H$ of unknown length, for elements that 
    are at distance at most  $r_{BCP}(\lambda)$. If $M_H$ is the bound of the area for the thin combing of $H$, 
    then we see that the area of $\check(D_i)$ is at most  
    $A_{\mathcal{K}} (M\times (1+  r_{BCP}(\lambda))) + 2 r_{BCP}(\lambda)  M_H$. 
    This proves the lemma, since there are at most $M$ such loops.

  \end{proof*}

  We now prove that $\check{\rho}$ is a  
  straight thin combing. That it is a thin combing follows 
  from the previous lemma and Lemma \ref{lem;ext_qg}. The quasi-geodesy follows from Lemma \ref{lem;qg2}.  
  We prove now that it is almost prefix-closed. 

  Let $v$ be a vertex on $\mathcal{K}$, and $w$ be on $\check{\rho} (v)$. 
  Let $w^-$ be the last vertex of $\rho(v)$ that is also on $\check{\rho}(v)$. 
  The vertices $w$ and $w^-$ are neighbors of a same coset vertex in $\widehat{\mathcal{K}}$. Therefore, 
  they are at distance at most $2$ in this graph, and we deduce that 
  $Area_{\widehat{\mathcal{K}}}(\rho(w),\rho(w^-),[w^-,w] ) \leq 2M$. Moreover, by prefix-closed property for $\rho$, one has 
  $Area_{\widehat{\mathcal{K}}}(\rho(w^-),\rho(v)|_{[v_0,w^- ]} ) \leq M$. Therefore, 
  $Area_{\widehat{\mathcal{K}}}( \rho(v)|_{[v_0,w^- ]}, [w^-,w], \rho(w)  ) \leq 3M$. Now notice that, by definition of $w^-$, 
  the extension of the concatenation $\rho(v)|_{[v_0,w^- ]}\cdot [w^-,w]$ is precisely $\check{\rho}(v)|_{[v_0,w]}$. 
  The previous lemma concludes the proof.

\end{proof}

\section{Free constructions.}

  In this section we explain how to construct straigth thin combings for 
  free constructions of groups.

\begin{theorem}\label{theo;free}
  
    Let $Y$ be a finite graph of groups, and $T$ the Bass-Serre tree, and $\Gamma$ be $\Gamma= \pi_1(Y,*)$. 
 Let $G_1 \dots G_N$ be the groups of vertices of $Y$. 
 
  (i)  We assume that, for each $G_i$, there exists a  $\Gamma$-graph of finite type $X^i$, with a $G_i$-$\lambda$-straight thin combings of $X^i$

  (ii) We assume also that each edge group of $Y$ fixes one and only one 
   vertex in each of these graphs in which it acts.

  Then, there is a $\Gamma$-graph of finite type that admits a $\Gamma$-$\lambda$-straight thin combing, and whose  
  stabilisers of the vertices of infinite valence are the conjugates 
  of the peripheral subgroups of the vertex groups. 

\end{theorem}

\begin{proof}

  We begin by constructing a graph on which the group $\Gamma$ acts. It is in fact a tree of spaces, constructed with 
  the Rips polyhedra $P(X^i)$  of the graphs $X^i$, in the way of the Bass-Serre tree, but we now define it precisely.  
  We denote by $\mathcal{K}^i$ the 1-skeleton of $P(X^i)$.
  For all $v$ of $T$ we choose an element $\gamma_v \in \Gamma$ such that $\gamma_v^{-1} v \in \{v_0, \dots v_N \}$.

Let $\tilde{P}$ be the complex  
  $\tilde{P} = \bigcup_{v\in T^0} \{ \gamma_v \} \times P(X^{i(v)})$, where for all $v$, $i(v)$ 
  is such that $\gamma_v^{-1} v = v_{i(v)}$
 The group $\Gamma$ acts on the complex $\tilde{P}$ as follows: 
  $\gamma (  \gamma_v,  p) = (\gamma_{\gamma v}, \gamma \gamma_v  (\gamma_{\gamma v})^{-1} p)$, 
  as can be easily checked.

  For all edge $\epsilon=(v,v')$ of $T$, the stabiliser $G_e$ of $e$ fixes one and only one vertex in 
  both $\{ \gamma_v\} \times P(X^{i(v)}) $ and  $\{ \gamma_v'\} \times P(X^{i(v'))}$. 
  We identify these two points, and we do that for each edge $\epsilon$. 
  Let $\mathcal{K}$ be the 1-skeleton of the complex obtained this way, $\mathcal{K} = (\tilde{P} / \sim)^{(1)}$, 
  where $\sim$ is the 
  identification mentionned. It is easily seen that it is a connected graph, and, since $T$ is a 
  tree, that for all vertex $v$ of $T$, the subgraph $\{ \gamma_v\} \times \mathcal{K}^{i(v)} $ of $\tilde{K}$ 
  embeds injectively in $\mathcal{K}$ in the quotient by $\sim$. We still denotes its image in $\mathcal{K}$ by 
  $\{ \gamma_v\} \times \mathcal{K}^{i(v)} $. 

  It is easily checked that the graph $\mathcal{K}$ contains finitely many orbits of edges, and that every edge has a vertex with finite stabiliser.  
 The vertices of infinite valence are the translates of the vertices of infinite valence of the graphs $\mathcal{K}^i$.

  Let $p_0$ be the point $(v_0, p_0^0)$ of $\mathcal{K}$. 

  For all vertex $p$ in $\mathcal{K}$. We define 
  $\mathcal{V} (p) = \{ v \in T, p\in  \{ \gamma_v\} \times \mathcal{K}^{i(v)}   \}$, and 
  $v(p)$ to be the closest to $v_0$ in this set. As $T$ is a tree, $v(p)$ is well defined and unique, for all $p$.

    \begin{lemma}
      For all $v\neq v_0$ there exists a vertex $p_e(v)$ in $\{ \gamma_v\} \times \mathcal{K}^{i(v)}$, such that 
      for all $p\in \{ \gamma_v\} \times \mathcal{K}^{i(v)}$, any path from $p_0$ to $p$ contains $p_e(v)$.
    \end{lemma}

    \begin{proof*}[Proof.]
 
      Let $v^-$ the vertex of $T$ that precede $v$ on the segment $[v_0,v]$ of $T$, and let  
      $p_e(v)$ be the unique vertex of $\mathcal{K}$ that is fixed by the stabilizer of the edge $(v^-,v')$.  
      By construction of $\mathcal{K}$, $p_e(v)$ is the only common point of $\{\gamma_{v^-} \}\times \mathcal{K}^{i(v^-)}$ and  
      $\{\gamma_{v} \}\times \mathcal{K}^{i(v)}$. Moreover, any path in $\mathcal{K}$ projects on a 
      path in $T$ by the following  
      rule: the image of any vertex $p'$ not fixed by an edge is  
      the unique vertex of $\mathcal{V}(p')$, and the image of a vertex $p'$ fixed by an edge stabilizer is the segment  
      $[ v(p'^-), v(p'^+)] \subset \mathcal{V}(p')$ of $T$, where $p'^-$ and $p'^-$ are the vertex of the path  
      preceding and following $p'$ (note that they are not fixed by an edge stabiliser), and   
      where $v(p'^-)$ is the unique vertex of $T$ in 
      $\mathcal{V}(p'^-)$ (similarily for $v(p'^+)$). As $T$ is a tree, we see that any path from $p_0$ to a vertex of    
      $\{ \gamma_v\} \times \mathcal{K}^{i(v)}$ contains a vertex of  
      $\{ \gamma_{v^-}\} \times \mathcal{K}^{i(v^-)}$, and therefore contains $p_e(v)$. 

    \end{proof*}

  We now define $\rho(p)$ for all $p \in \mathcal{K}$. 
  Consider the segment $[v_0,v(p)]$ in $T$, and denote by $w_0=v_0, w_1, \dots w_r, w_{r+1} = v(p)$ the 
  consecutive vertices of this segment in $T$. Then, 
  $$ \rho( p) = \rho_0^0 (p_e(w_1) ) \, \cdot \, \gamma_1 \rho_{i_1}^{j_1}(\gamma_1^{-1} p_e(w_2)) 
  \,\cdot  \dots  \cdot \, \gamma_{r-1} \rho_{i_{r-1}}^{j_{r-1}}(\gamma_{r-1}^{-1} p_e(w_r)) \, \cdot  
  \gamma_r \rho_{i_r}^{j_r} (\gamma_r^{-1} p) $$
  where $\cdot$ means concatenation, and where for all $k$, $\gamma_k$, $i_k$ and $j_k$ are the only possible 
  values so that the formula makes sense, that is so that $\gamma_k^{-1} w_k = v_{j_k} \in \{ v_0, \dots, v_N \}$, 
  and   $\gamma_k^{-1} p_e(w_k) = p^{j_k}_{i_k} \in \mathcal{K}^{j_k}$.

  We check now that $\rho$ is a thin combing. Let $p$ and $p'$ be neighbors in $\mathcal{K}$.

  Let us first assume that $v(p)=v(p')$. 
  By definition of $\rho$, we see that $Area (\rho(p), (p,p'), \rho(p'))$ is equal to 
  $Area ( \gamma_r \rho_{i_r}^{j_r} (\gamma_r^{-1} p), (p,p'),  \gamma_r \rho_{i_r}^{j_r} (\gamma_r^{-1} p))$, 
  because the paths $\rho(p)$ and $\rho(p')$ coincide until $\gamma_r p_{i_r}^{j_r}$. 
  But this area is bounded by $M_{i_r}^{j_r}$, the constant for the thin combing  $\rho_{i_r}^{j_r}$.

  If neither $p$ nor $p'$ is fixed by the stabilisor of an edge of $T$, 
  then it is clear that $v(p)=v(p')$, and we have treated this case.  

   The vertices  $p$ and $p'$ cannot be both fixed by the stabilisor of an edge of $T$, since such two points 
   are at distance at least $2$, by assumption on the graphs $\mathcal{K}^i$. 
  If now $p$ is fixed by the stabilisor of an edge of $T$, and not $p'$, and if $v(p)\neq v(p')$ 
  (they need not be neighbors in $T$). Let us assume that $v(p)$ is closer to $v_0$ than $v(p')$. As $p$ and $p'$ 
  are neighbors, $\{ \gamma_v(p)\} \times \mathcal{K}^{i(v(p))}$ and 
  $\{ \gamma_v(p')\} \times \mathcal{K}^{i(v(p'))}$  share one point, and this is necessarily $v(p)=p$. 
  By definition, $\rho(p') =  \rho_0^0 (p_e(w_1) ) \, \cdot \, \dots , \cdot \gamma_r \rho_{i_r}^{j_r}(\gamma_r^{-1} p_e(w_r)) 
  \,\cdot \, \dots \, \cdot  \gamma_s \rho_{i_s}^{j_s} (\gamma_s^{-1} p')$, but we see that the paths   
  $\gamma_{k} \rho_{i_{k}}^{j_{k}}(\gamma_{k}^{-1} p_e(w_{k}))$, for $k=r\dots s-1$, are all paths from $v_p$ to $v_p$, 
  being quasi-geodesics, they are therefore all trivial. Hence, $Area (\rho(p), (p,p'), \rho(p'))$ is equal to
    $Area ((p,p') \gamma_s \rho_{i_s}^{j_s} (\gamma_s^{-1} p'))$, which is bounded above by property of  
  $\rho_{i_s}^{j_s}$. Hence, $\rho$ is a thin combing.

  Let us prove now that it is almost prefix-closed. Let $p$ be a vertex of $\mathcal{K}$ and $p'$ be a vertex of $\rho(p)$. 
  Let $\rho(p)=  \rho_0^0 (p_e(w_1) ) \, \cdot \, \gamma_1 \rho_{i_1}^{j_1}(\gamma_1^{-1} p_e(w_2)) 
  \,\cdot \, \dots \, \cdot  \gamma_r \rho_{i_r}^{j_r} (\gamma_r^{-1} p) $, by definition.  
 
  If the point $p'$ is in $\gamma_r \rho_{i_r}^{j_r} (\gamma_r^{-1} p)$, and is not the first point $\gamma_r p_{i_r}^{j_r}$ 
  of this  
  subpath, then $\rho(p')$ has same prefix than $\rho(p)$ until  $\gamma_r p_{i_r}^{j_r}$, and therefore,  
  $Area (\rho(p)|_{[p_0,p']}, \rho(p'))$ is 
  equal to  $Area (\gamma_r \rho_{i_r}^{j_r} (\gamma_r^{-1} (w_{r+1})|_{[p_{i_r}^{j_r} ,p']} , 
  \gamma_r \rho_{i_r}^{j_r} (\gamma_s^{-1} p')))$, 
  which is bounded above by property of prefix-closedness for $\rho_{i_r}^{j_r}$. 
 
  If now $p'$ is not as before, then there exists a unique index $s<r$ such that  
  $p' \in \gamma_s \rho_{i_s}^{j_s}(\gamma_s^{-1} p_e(w_{s+1}))$ and $p' \neq \gamma_s p_{i_s}^{j_s}$.   
  Then, $Area (\rho(p)|_{[p_0,p']}, \rho(p'))$ is equal to  
  $Area (\gamma_s \rho_{i_s}^{j_s} (\gamma_s^{-1} 
  (w_{s+1})|_{[p_{i_s}^{j_s} ,p']} , \gamma_s \rho_{i_s}^{j_s} (\gamma_s^{-1} p')))$, 
  which is bounded above by $M$ by property of prefix-closedness for $\rho_{i_s}^{j_s}$. This proves that $\rho$ 
  is almost prefix closed.

  It remains to see that for all $p$, the path  $\rho(p)$ is a $\lambda$ quasi-geodesic. Again it is an application 
  of the lemma, and of the quasi-geodicity of the pieces of $\rho(p)$.

\end{proof}

\section{Applications to limit groups.}

We now give applications of the sections 2 and 3 to the case of limit groups,
and of groups acting freely on $\mathbb{R}^n$-trees. 

First we note that the work of the previous part allows to consider some free constructions over cyclic groups that are not peripheral.
Before stating the corallary, let us recall

\begin{coro} \label{coro;comb_f}

    Let $Y$ be a finite graph of groups, and $T$ the Bass-Serre tree. 
  Let us choose a set of representatives of orbits of vertices in $T$ by the action of $\Gamma= \pi_1(Y,*)$, 
  for a base point $*$ in $Y$: $\{v_0, \dots v_N \} $. 
  For all $i =0 \dots N$, let $G_i$ be the stabiliser of $v_i$ in $\Gamma$. 
 
   We assume that, for each $G_i$, there is a hyperbolic and fine $G_i$-graph of finite type, $\mathcal{K}^i$.  

  We assume that every edge group is either finite or cyclic. If an edge group is cyclic, we assume that it is maximal cyclic in at least one of the vertex groups of the edge.

  If, there exists, for every vertex group $G_i$,  a $G_i$-straight thin combing of $\mathcal{K}^i$, then there exists a $_Gamma$-graph of finite type, such that the stabilisers of vertices of infinite valence are the conjugates of the peripheral groups of vertex groups, and that admits a $\Gamma$-straight thin combing.

\end{coro}

\begin{proof}

Remark. By the Combination Theorem of \cite{Dgt}, the group $\Gamma$ is hyperbolic relative to the conjugates of the peripheral subgroups of the vertices groups. Let $\mathcal{G}$ be this family of peripheral subgroups. Let $\mathcal{H}$ be the family of the stabilisers of edges of the Bass Serre tree that are not in $\mathcal{G}$.

 Let $v_i$ be a vertex of the graph of groups $X$. 
 Let us consider the conned off graph of $\mathcal{K}^i$  along the orbits of every edge group adjacent to $v_i$, 
 and denote it by $\widehat{\mathcal{K}^i}$. By Lemma \ref{lem;contrac}, there exists a $G_i$-straight thin-combing of $\widehat{\mathcal{K}^i}$. By Theorem \ref{theo;free}, there  exists a $\Gamma$-straight-thin combing on a $\Gamma$-graph of finite type $\mathcal{K}$,  
 whose stabilisers of infinite vertices are the elements of $\mathcal{G}\cup \mathcal{H}$. 
 Let $P$ be a polyhedron in which the areas of this thin combing are measured. As it is locally finite around the edges, it is contained in any complex $P_{d,r}(\mathcal{K})$, for $d$ and $r$ large enough.  
 By relative hyperbolicity of $\Gamma$, relatively to $\mathcal{G}$,  and by a Theorem in \cite{Dplms} (Theorem 6.2),
there exists a relative Rips polyhedron for $\check{\mathcal{K}}$, in particular, it is 1-connected. 
 Such a complex is contained in $P_{d,r}(\mathcal{K})$, if $d$ and $r$ are large enough.
 Therefore, without loss of generality, we can assume that $P$ contains the 2-coned-off of a Rips polyhedron of  $(\Gamma, \mathcal{G})$,   
such a polyhedron satisfies an isoperimetric inequality, because it is simply connected and has finitely many orbits of simple loops of given length (see Theorem 6.2 \cite{Dplms}). 

Therefore we can apply Lemma \ref{lem;ext} to get that $\Gamma$ acts on a $\Gamma$-graph of finite type whose stabilisers of infinite vertices are the elements of $\mathcal{G}$, and that admits a $\Gamma$-straight-thin combing.
 
\end{proof}

\begin{coro}

 If $\Gamma$ acts freely on an $\mathbb{R}^n$-tree (in particular if it is a limit group), 
  then it acts on a $\Gamma$-graph of finite type whose stabilisers of infinite vertices are the maximal abelian non cyclic subgroups, and that admits a $\Gamma$-straight thin combing.

\end{coro}

\begin{proof}

This is a consequence of the D\'evissage Theorem of Guirardel (Theorem 7.1 \cite{VG}, or Coro. 6.6 in \cite{Gros}, see also Sela \cite{S} for limit groups)  and of the previous corollary \ref{coro;comb_f}. Indeed, the D\'evissage theorem and the Corollary \ref{coro;comb_f}, allow us to reduce the problem to groups acting freely on $\mathbb{R}^{n-1}$-trees, and by induction, it is enough to solve the problem for groups acting on $\mathbb{R}$-trees. By a theorem of Rips,  these groups are abelian or surface groups. and therefore, by the corollary  \ref{coro;comb_f}, it is enough to consider only abelian groups and free groups (or even cyclic groups), and the result is then obvious.

\end{proof}

\begin{coro}

  Any limit group is metabolic relative to its maximal non-cyclic abelian subgroups.

\end{coro}

\begin{proof}

 We have to show that any 2-$\ell_{\infty}$-cocycle of a relative Rips polyhedron with values in a normed abelian group 
 is the co-boundary of a 1-$\ell_{\infty}$-cocycle. Let $c$ be such a 2-cocycle. We define $c'$ on the edges of the Relative Rips polyhedron to be $c'((v_1,v_2)) = \sum_{\sigma \in (\rho(v_1), \rho(v_2), (v_1,v_2))} c(\sigma)$, where $\rho$ is a thin combing, and where the sum is performed on a coice of 2-cell that fills the loop $(\rho(v_1), \rho(v_2), (v_1,v_2))$ in away that minimizes area. By definition of the thin combing, and by boundedness of $c$, $c'$ is a $\ell_{\infty}$-cocycle. For every 2-cell $f$,  $\partial c'(f)$ is equal to $c(f)$ by cocycle property for $c$. This means that $\partial c'=c$, what we needed. 

\end{proof}

\begin{coro}

  Any hyperbolic limit group is metabolic.

\end{coro}

Such a group acts cocompactly on a locally finite graph that admits a thin combing, this implies metabolicity from \cite{G2}, or from the above.

\end{document}